\documentclass[11pt]{amsart}
\usepackage{graphics}
\usepackage{eepic}
\usepackage{amsmath}
\usepackage{amscd}
\usepackage{amsfonts,latexsym,amssymb}

\newcommand{\SA}{{\mathcal{A}}}

\newcommand{\SD}{{\mathcal{D}}}
\newcommand{\SE}{{\mathcal{E}}}
\newcommand{\SF}{{\mathcal{F}}}

\newcommand{\SH}{{\mathcal{H}}}

\newcommand{\SO}{{\mathcal{O}}}
\newcommand{\SP}{{\mathcal{P}}}

\newcommand{\SR}{{\mathcal{R}}}

\newcommand{\CrD}{{\SD}}
\newcommand{\CrE}{{\SE}}
\newcommand{\CrF}{{\SF}}
\newcommand{\CrH}{{\SH}}

\newcommand{\CrO}{{\SO}}
\newcommand{\CrP}{{\SP}}
\newcommand{\CrR}{{\SR}}

\newcommand{\AAA}{\mathbb{A}}
\newcommand{\HH}{\mathbb{H}}

\newcommand{\CC}{\mathbb{C}}

\newcommand{\isom}{\cong}
\newcommand{\Ext}{\operatorname{Ext}}
\newcommand{\SEnd}{End}

\newcommand{\SHom}{Hom}

\newcommand{\id}{\textbf{\textit{I}}}
\newcommand{\im}{\operatorname{im}}

\newcommand{\rk}{\operatorname{rk}}

\newcommand{\too}{\longrightarrow}

\newcommand{\gln}{{Gl_n(\CC)}}
\newcommand{\glp}{{Gl_{n+1}(\CC)}}
\newcommand{\afn}{{Af_{n}(\CC)}}
\newcommand{\sln}{{SL_n(\CC)}}
\newcommand{\matn}{{Mat}_{n \times 1}}
\newcommand{\un}{{U(n)}}
\newcommand{\holoextn}{\underline{{\CrE}xtn}}
\newcommand{\extn}{\underline{Extn}}
\newcommand{\flatextn}{\underline{FlatExtn}}
\newcommand{\higgsextn}{\underline{{\CrH}iggs{\CrE}xtn}}
\newcommand{\pairs}{\underline{{\CrP}airs}}
\newcommand{\flatp}{\underline{\CrD\CrR\CrH}}
\newcommand{\thw}{\theta_\wedge}
\newcommand{\delb}{{\overline\partial}}
\newcommand{\uni}{\widetilde {\underline U}}

\newtheorem{proposition}{Proposition}[section]
\newtheorem{theorem}[proposition]{Theorem}
\newtheorem{definition}[proposition]{Definition}
\newtheorem{lemma}[proposition]{Lemma}

\newtheorem{remark}[proposition]{Remark}

\author[T. G\'omez and F. Presas]{Tom\'as L. G\'omez 
and Francisco Presas }
\title[Affine representations]{Affine
representations of the fundamental group
(with an appendix on parabolic representations)}
\thanks{T.G. was supported by a postdoctoral fellowship of
Ministerio de Educaci\'on y Cultura (Spain). F.P. was supported by a
postgraduate fellowship of Universidad Complutense de Madrid.}


\date{11 November 1999}

\address{Tom\'as L. G\'omez:
School of Mathematics,
Tata Institute of Fundamental Research,
Homi Bhabha Road, 400 005 Mumbai (India)
}

\email{tomas@math.tifr.res.in}

\address{Francisco Presas:
Departamento de Algebra, 
Facultad de Ciencias Ma\-te\-m\'a\-ti\-cas,
Universidad Complutense de Madrid, 28040 Madrid (Spain)
}

\email{fpm@eucmos.sim.ucm.es}

\begin{document}

\begin{abstract}
Unitary representations of the fundamental group of a K\"ahl\-er manifold
correspond to 
polystable vector bundles (with vanishing Chern classes). 
Semisimple linear representations correspond to polystable Higgs
bundles. In this paper we find the objects corresponding to affine
representations: the linear part gives a Higgs bundle and
the translation part corresponds to an element of a generalized de
Rham cohomology.
\end{abstract}

\maketitle

\section{Introduction}

The affine group is the subgroup of $\glp$
$$
\afn= \left\{ 
\left(
\begin{array}{cc}
A & t \\
0 & 1
\end{array}
\right) 
: \quad A \in \gln,\quad t\in \matn \right\}.
$$
Given an element of $\afn$, the matrix $A\in \gln$ is called the
linear part and the column vector $t\in \matn$ is called the
translation part. The map that gives the linear part defines a group
homomorphism $\afn \to \gln$.

Recall that given a K\"ahler manifold $X$, the set of conjugacy classes of
unitary representations of the fundamental group 
$\pi_1 \to \un$ is equal to the set of
isomorphism classes of polystable rank $n$ vector bundles with
vanishing Chern classes (\cite{NS}, \cite{D}, \cite{UY}). 
This correspondence can be extended to
semisimple $\gln$ representations. Then we have to consider polystable
rank $n$ Higgs bundles (\cite{H}, \cite{S1}, \cite{S2}, \cite{S3}). 
In this paper we study affine representations.

Given an affine representation, the linear part gives a $\gln$
representation, and if this is semisimple, this defines a polystable
Higgs bundle $(\CrE,\theta)$. In this paper we prove that the extra
object we have to
add to the Higgs bundle to account for the translation part is an
element of the first de Rham cohomology group $H^1_{DR}(E)$. This
cohomology was introduced by Simpson in \cite{S1}.

{}From a different point of view, since $\afn \subset \glp$, an affine
representation gives a $\glp$ representation, but unless the
translation part $b$ is equal to zero, this representation won't be
semisimple. Arbitrary (not necessarily semisimple) representations have
been studied by Simpson in \cite{S1} using differential graded
categories. Roughly speaking, he shows that it is necessary 
to add a new field
$\eta \in A^1(X,\SEnd(E))$, a smooth 1-form with values in
$\SEnd(E)$, where $E$ is the smooth vector bundle underlying the 
Higgs bundle 
corresponding to the induced
semisimple representation.

{}From this point of view, what we show in this paper is that in the
case of affine representations (assuming that the $\gln$ representation
given by the linear part $A$ is semisimple) there is an explicit
cohomological interpretation of $\eta$ in terms of de Rham
cohomology. If furthermore the linear part $A$ is unitary, de Rham
cohomology can be described in terms of the usual cohomology groups of
coherent sheaves.

The parabolic construction of principal
$G$-bundles on an elliptic curve done by \cite{FMW} is related to
this. Given a certain maximal parabolic subgroup $P$ of $G$, 
with Levi factor
$L$ and unipotent part $U$, first they construct a principal
$L$-bundle (the semisimple part of $P$), and then they show that the
extra piece of data needed to specify the $P$-bundle is an element in
the \'etale cohomology group $H^1_{et}(X,\underline U)$, where
$\underline U$ is the associated principal $U$-bundle. In the case $G=
\sln$, this is a usual extension group $\Ext^1$ of
vector bundles.

The techniques used to study affine representations of the fundamental
group can be adapted 
to consider representations into a parabolic subgroup of $\gln$. 
The details are
given in an appendix.

\noindent\textbf{Notation.} Holomorphic vector bundles will be denoted 
$\CrE$, $\CrF$, ... and the corresponding underlying smooth vector
bundles will be denoted $E$, $F$, ... The trivial holomorphic line
bundle is denoted $\CrO_X$, and the underlying smooth bundle $O_X$.

\bigskip

\noindent\textbf{Acknowledgments.} We would like to thank Ignacio
Sols, for proposing us to study affine representations of the
fundamental group, and for
several discussions. We would also like to thank M.S. Narasimhan for
discussions on this subject. T.G. thanks Universidad Complutense de
Madrid for its hospitality during the visit where this work was 
started.

\section{Affine bundles}
\label{affinebundles}

A holomorphic (resp. smooth) affine bundle is a holomorphic
(resp. smooth) principal bundle 
with structure group
$\afn$. Given a principal $\afn$-bundle, by the inclusion $\afn
\subset \glp$ we obtain a principal $\glp$-bundle. Since the
transition functions $\{\alpha_{ij}\}$ of this principal bundle are of
the form
\begin{eqnarray}
\label{transfun}
\alpha_{ij}=
\left(
\begin{array}{cc}
A_{ij} & t_{ij} \\
0 & 1
\end{array}
\right), 
\end{eqnarray}
the associated rank $n+1$ vector bundle 
$\CrF$ has a 
canonical rank $n$
subbundle $\CrE$ (whose transition functions are $\{A_{ij}\}$), and the
quotient is isomorphic to the trivial line bundle $\CrO_X$ (since the
right lower entry of $\alpha_{ij}$ is equal to 1)
$$
0 \to \CrE \to \CrF \to \CrO_X \to 0.
$$
Conversely, given an extension like this, we can find local
trivializations of $\CrF$ such that the transition functions are of the
form (\ref{transfun}), and hence we obtain a principal
$\afn$-bundle. It is easy to check that the isomorphism class of this
affine bundle doesn't
depend on the choice of local trivializations.

\begin{definition}
Let $\holoextn$ (resp. $\extn$) be the category
 whose 
objects are 
holomorphic (resp. smooth) extensions
of the form 
\begin{eqnarray}
\label{obj}
0 \to \CrE \to \CrF \to \CrO_X \to 0
\end{eqnarray}
(with $\rk(\CrF)=n+1$, $\rk(\CrE)=n$), and whose morphisms are pairs
$(\varphi,\psi)$ where $\varphi:\CrE\to \CrE'$ and $\psi:\CrF \to \CrF'$ are
vector bundle morphisms and the following diagram commutes
\begin{eqnarray}
\label{morph}
\CD
0 @>>>  \CrE            @>{i}>> \CrF  @>{p}>>  \CrO_X @>>> 0 \\
@.     @V{\varphi}VV         @V{\psi}VV   @|     @. \\
0 @>>>  \CrE'          @>{i'}>> \CrF' @>{p'}>> \CrO_X @>>> 0. \\
\endCD
\end{eqnarray}
\end{definition}

\begin{lemma}
\label{affinelemma}
The category of
holomorphic (resp. smooth) principal $\afn$-bundles 
is equivalent to the category $\holoextn$ (resp. $\extn$).
\end{lemma}

\hfill $\Box$

The proof is analogous to the proof of the fact that the
category of principal $\gln$-bundles is equivalent to the category of
rank $n$ vector bundles.

\begin{remark}
\textup{
One could be tempted to relax condition (\ref{morph}), and to allow
for an arbitrary isomorphism on $\CrO_X$, instead of requiring it to be
the identity (as it is done for the definition of weak isomorphism of
extensions). But this wouldn't give the right category. This is easily
checked by looking at the case where the base space $X$ is a
point. The dimension of the automorphisms group should then be
$n^2+n=\dim (\afn)$. This is the dimension that we obtain with the
definition that we have given, but if we allowed for an arbitrary
isomorphism on $\CrO_X$, the dimension would be $n^2+n+1$.}

\textup{
But if we are only interested in the set of isomorphism classes, this
distinction is not important, since any isomorphism of $\CrO_X$ is
multiplication by scalar, and this can be absorbed by rescaling
$\varphi$ and $\psi$.}
\end{remark}

\begin{remark} \label{affinemodel}
\textup{
Abstractly, an affine space modeled on a vector space $V$ is a set
$\AAA$ together with a transitive and free action of $V$, and the
affine group is the automorphism group of $\AAA$.
Given an extension as in (\ref{obj}), note that $\CrE \isom p^{-1}(0)
\subset \CrF$ is a bundle of vector spaces that acts (by addition) on
$A=p^{-1}(1)\subset \CrF$
$$
\CrE \times A \to A,
$$
and the action commutes with projection to the base space $X$.
Then we can think of $A$ as a bundle of affine spaces. Conversely,
given $A$ and a vector bundle $\CrE$ acting on $A$, we can recover (up to
isomorphism) the extension $\CrF$. This gives an equivalent description
of an affine bundle.}
\end{remark}

\begin{definition}
Let $\pairs$ be the category whose objects are
pairs $(\CrE,\chi)$ where $\CrE$ is a rank $n$ vector bundle and $\chi\in
H^1(\CrE)$, and a morphism between $(\CrE,\chi)$ and $(\CrE',\chi')$ is a
morphism of vector bundles $\varphi:\CrE\to \CrE'$ such that 
$$
\begin{array}{ccc}
H^1(\CrE) & \stackrel{H^1(\varphi)}\too & H^1(\CrE') \\
\chi &\longmapsto &\chi'\quad .
\end{array}
$$
\end{definition}

An object (\ref{obj}) of $\holoextn$ gives an element $\chi \in
\Ext^1(\CrO_X,\CrE)=H^1(\CrE)$, and hence an element of $\pairs$.
Note that if $\lambda\neq 0$ then $(\CrE,\chi)$ and $(\CrE,\lambda\chi)$ are
isomorphic (take $\varphi=\lambda \id_\CrE$). This category
is not equivalent to the category of affine bundles, but we have

\begin{proposition}
There is a natural bijection between the isomorphism classes of
principal holomorphic $\afn$-bundles and the isomorphism classes 
of pairs $(\CrE,\chi)$,
where $\CrE$ is a rank $n$ vector bundle and $\chi\in H^1(\CrE)$.
\end{proposition}

\begin{proof}
Using lemma \ref{affinelemma}, the set of isomorphism classes of
holomorphic principal
$\afn$-bundles is equal to the set of isomorphism classes of the
category $\holoextn$. Two isomorphic extensions as in (\ref{morph}) give
isomorphic pairs $(\CrE,\chi)$ and  $(\CrE',\chi')$, 
because (\ref{morph}) implies that the image of $\chi$ in $H^1(\CrE')$
under the map induced by $\varphi$ is $\chi'$.

Conversely, two isomorphic pairs give two extensions (unique up to
noncanonical isomorphism) that are isomorphic in the category
$\holoextn$.

\end{proof}

\begin{remark}
\textup{This proposition is also valid in the smooth category. 
Note that in the smooth category short exact
sequences are always split, although this splitting is not
unique
(in the smooth category $H^1(E)$ is
zero since the sheaf of smooth sections of $E$ is fine). 
This means that a smooth affine
bundle has a (smooth) reduction of structure group to $\gln
\subset \afn$, where this inclusion is given by
$$
\left(
\begin{array}{cc}
A & 0 \\
0 & 1 
\end{array}
\right)
\subset
\left(
\begin{array}{cc}
A & t \\
0 & 1 
\end{array}
\right)
$$
The existence of this reduction is equivalent to the fact that any
affine bundle has a smooth section. From a topological point of view if
we interpret
an affine bundle, following remark \ref{affinemodel}, as a fibration of
affine spaces
over $X$ the existence of the smooth section is a consequence of the
contractibility of
the fibers.
}
\end{remark}

\section{Flat affine bundles}
\label{flataffinebundles}

A (smooth) flat affine bundle is a (smooth) principal $\afn$-bundle 
$A$ with a connection $D^{}_A$ such that 
$D^2_A=0$. A morphism of flat
affine bundles is a morphism $f:A\to A'$ of principal bundles such
that the pullback of the connection on $A'$ is equal to the connection
on $A$.

\begin{definition} 
Let $\flatextn$ be the category 
whose objects are extensions
\begin{eqnarray}
0 \to E \to F \to O_X \to 0
\label{extension}
\end{eqnarray}
with $\rk(E)=n$ and $\rk(F)=n+1$,
together
with a flat connection $D_F$ on $F$ that respects $E$ (in the sense
that the image of $\SA^0(X,E)$ is in
$\SA^1(X,E)\subset \SA^1(X,F)$, and then it induces a connections on
$E$ and $O_X$) and that induces
the trivial connection $D_{O_X}$ on $O_X$. A morphism in
this category is a morphism in $\extn$ with $\varphi^*D_{F'}=D_F$.
\end{definition}

An affine bundle gives an extension like (\ref{extension}) by lemma
\ref{affinelemma}. A flat affine connection gives a connection
$$
D_F: \SA^0(X,F) \to \SA^1(X,F)
$$
that preserves $E$ , and the induced connection on $O_X$ is the trivial
connection $D_{O_X}$, and then we get an object of $\flatextn$. In
fact, this construction gives an equivalence of categories:

\begin{lemma}
\label{flataffinelemma}
The category of flat affine bundles is equivalent to the category
$\flatextn$.
\end{lemma}

\hfill $\Box$

If the induced flat connection $D_E$ on $E$ is semisimple, then by 
\cite[Theorem 1]{S1} there
is a harmonic metric on $E$ and a decomposition $D_E=\nabla_E + \alpha$
where $\nabla_E$ is a unitary connection and $\alpha\in A^1(X,\SEnd(E))$ 
is a 1-form. Let $\partial^{}_E$ and $\delb^{}_E$ be the (1,0) and
(0,1) part of
$\nabla^{}_E$. We get a holomorphic
vector bundle
$\CrE=(E,\delb^{}_E)$. Since the metric is harmonic
the 1-form can be written as $\alpha=\theta +\theta^*$ where $\theta \in
H^0(\SEnd(\CrE)\otimes \Omega^1)$ is a holomorphic (1,0) form 
with values in
$\SEnd(\CrE)$ such that $\theta_\wedge \theta=0$, and $\theta^*$ is the
conjugate (0,1) form. Then $\theta$ is a Higgs field, and
the Higgs bundle $(\CrE,\theta)$ is polystable.
Following Simpson, we decompose $D_E'=\partial^{}_E + \theta^*$, 
$D_E''=\delb^{}_E + \theta$, and define the following cohomology groups:

\begin{itemize}

\item The de Rham cohomology $H^i_{DR}(E)$ of a Higgs bundle is the 
cohomology
of the complex $(A^\bullet(X,E),D)$, where $A^p(X,E)$ is the space
of global p-forms with coefficients in $E$.

\item The Dolbeault cohomology $H^i_{Dol}(E)$ is the
cohomology of the complex $(A^\bullet(X,E),D''_E)$.

\item The group $\HH^i(E\otimes \Omega^\bullet)$ is defined as the 
hypercohomology of
the complex
\begin{eqnarray}
\label{holomcomplex}
E \stackrel{\thw}\longrightarrow E\otimes \Omega 
\stackrel{\thw}\longrightarrow E\otimes \Omega^2 
\stackrel{\thw}\longrightarrow \cdots
\end{eqnarray}
\end{itemize}
Simpson shows in \cite{S1} that these three cohomologies are naturally
isomorphic.

\begin{definition}
Let  $\higgsextn$ be the category whose objects are pairs
$$
(0 \to \CrE \to \CrF \to \CrO_X \to 0,\Theta)
$$
where the extension is an object of $\holoextn$ such that the Chern
classes of $\CrF$ vanish, and $\Theta\in H^0(\SEnd(\CrF)\otimes\Omega_X)$
is a Higgs field, such that $(\CrE,\Theta)$ is semistable, $\CrE$ is 
$\Theta$-invariant, and the Higgs field induced on $\CrO_X$ is zero.
A morphism in this category is a morphism in $\holoextn$ with 
$\psi^*\Theta'=\Theta$.
\end{definition}

\begin{definition}
Let $\flatp$ bet the category whose objects are triples
$(\CrE,\theta,\xi)$, where $(\CrE,\theta)$ is a polystable Higgs bundle with
vanishing Chern classes and $\xi \in H^1_{DR}(E)$. A morphism is a
morphism $\varphi$ of Higgs bundles such that
$$
\begin{array}{ccc}
H_{DR}^1(E) & \stackrel{H_{DR}^1(\varphi)}\longrightarrow & H_{DR}^1(E') \\
\xi &\longmapsto &\xi'\quad .
\end{array}
$$
\end{definition}

\begin{theorem}
\label{affinethm}
There are natural bijections among the following sets
\begin{enumerate}

\item
The set of $\afn$ representations of the fundamental group such that
the linear part $\gln$ is semisimple, modulo conjugation by elements
in $\afn$.

\item
The set of isomorphism classes of
objects of $\flatextn$ such that the induced flat connection $D_E$ is
semisimple
$$
\{ (0\to E \to F \to O_X \to 0,\; D_F): D_E\; \text{semisimple} 
\}/\isom\;.
$$

\item
The set of isomorphism classes of objects of $\higgsextn$ such that
the induced Higgs bundle $(\CrE,\theta)$ is polystable
$$
\{ (0\to \CrE \to \CrF \to \CrO_X \to 0,
\; \Theta): (\CrE,\theta)\; \text{polystable} \}/\isom\;.
$$

\item The set of isomorphism classes of $\flatp$,
$$
\{ (\CrE,\theta,\xi) \} /\isom\;.
$$
\end{enumerate}

\end{theorem}

\begin{proof}
$(1 \leftrightarrow 2)$ This follows from holonomy and
lemma \ref{flataffinelemma}.

$(2 \leftrightarrow 4)$ Take and object 
in $\flatextn$, i.e. an extension
$$
0 \to E \to F \stackrel{p}{\to} O_X \to 0
$$
with a flat connection $D_F$ preserving $E$. Take a $C^\infty$
splitting $E\oplus O_X \isom F$. This is given by a smooth 
morphism $\tau:O_X \to F$ (i.e., a smooth section of $F$) with $p\circ
\tau = \id_{O_X}$.
The fact that $D_F$ preserves $E$ and that it induces the
trivial connection $D_{O_X}$ on $O_X$ means that, using this
splitting, $D_F$ can be written as
$$
\left(
\begin{array}{cc}
D_E & b_\tau \\
0 & D_{O_X}
\end{array}
\right)
$$
where $D_E$ is the connection induced on $E$ and $b_\tau$ is a smooth
1-form with values in $E$. Note that $b_\tau$ depends on the splitting
(i.e. the smooth section $\tau$ chosen) but $D_E$ doesn't. Flatness of
$D_F$ translates into
$$
D_E^2=0 \quad \text{and} \quad D^{}_E b_\tau=0.
$$
Since by hypothesis 
$D_E$ is semisimple, as we have already explained it has a harmonic
metric and then there is a
polystable Higgs bundle $(\CrE,\theta)$ (this is
independent of the chosen section $\tau$, since $D_E$ was). If we take
a different splitting, i.e. change
the section $\tau$ to a new smooth section $\sigma$ of $F$ (again with
$p\circ \sigma=\id_{O_X}$), 
the diffeomorphism
between the old and the new splitting is given by a matrix of the form
$$
\left(
\begin{array}{cc}
\id_E & \tau-\sigma \\
0 & 1
\end{array}
\right)
$$
Note that the image of $\tau-\sigma$ is in $E$ (because $p\circ
(\tau-\sigma)=0$), and then this matrix makes sense. In the new 
splitting the connection $D_F$ has the form
$$
\left(
\begin{array}{cc}
D_E & b_\tau+D^{}_F (\sigma-\tau) \\
0 & D_{O_X}
\end{array}
\right)
$$
Then the 1-form associated to the new splitting is $b_\sigma=b_\tau
+D^{}_F (\sigma-\tau)$.
In other words, we obtain an element of $H^1_{DR}(E)$ (independent of
the splitting).

To show that this construction gives a bijection, we now construct an
inverse. Take a triple $(\CrE,\theta,\xi)$ where $(\CrE,\theta)$ is a
polystable Higgs bundle and $\xi\in H^1_{DR}(E)$. Since the Higgs
bundle is polystable, there is a harmonic metric that gives a
flat connection $D_E$.

Let $b\in A^1(X,E)$ be a representative of $\xi$ and define on 
$F=E \oplus O_X$ a connection
$$
D^{}_F=\left(
\begin{array}{cc}
D_E & b \\
0 & D_{O_X}
\end{array}
\right).
$$

$(4 \leftrightarrow 3)$ 
Let $(\CrE,\theta,\xi)$ be an object of $\flatp$. Let $\overline b
\in H^1_{Dol}(E)$ be the element corresponding to $\xi \in H^1_{DR}(E)$
under the natural isomorphism. Let $b \in A^1(X,E)$ be a $D''_E$-closed 
smooth 1-form
representing $\overline b$. Then 
\begin{eqnarray}
\label{closed}
\theta_\wedge b^{1,0}=0, \quad 
\delb_E b^{1,0} + \theta_\wedge b^{0,1}=0, \quad \text{and} \quad
\delb_E b^{0,1}=0.
\end{eqnarray}
Consider the vector bundle $F=E\oplus O_X$. Using this splitting,
let
$$
\delb^{}_F=
\left(
\begin{array}{cc}
\delb_E & b^{0,1} \\
0 & \delb_{O_X}
\end{array}
\right)
:A^0(X,F) \to A^{0,1}(X,F).
$$
Equations (\ref{closed}) imply that $\delb_F^2=0$,
and then this defines a structure of holomorphic vector bundle $\CrF$.
Let
$$
\Theta=
\left(
\begin{array}{cc}
\theta & b^{1,0} \\
0 & 0
\end{array}
\right)\in Hom(\CrF,\CrF\otimes\Omega).
$$
 Then
$\delb_{Hom(\CrF,\CrF\otimes\Omega)} \Theta=0$
and $\Theta_\wedge \Theta=0$ again by (\ref{closed}), and 
then $(\CrF,\Theta)$ is a 
Higgs bundle. It is easy to check that it is semistable (but not
polystable if $b^{1,0}\neq 0$).

Conversely, given an object of $\higgsextn$, let $F=E\oplus O_X$ be 
the underlying
smooth vector bundle of $\SF$, and let $\delb_F$ be the corresponding
$\delb$-operator. It can be written as
$$
\delb^{}_F=
\left(
\begin{array}{cc}
\delb_E & b_1 \\
0 & \delb_{O_X}
\end{array}
\right)
:A^0(X,F) \to A^{0,1}(X,F).
$$
with $b_1\in A^{0,1}(X,E)$ a smooth (0,1)-form. The fact that $\CrE$ is
$\Theta$-invariant and that it induces the zero Higgs bundle on
$\CrO_X$ imply that $\Theta$ can be written as
$$
\Theta=
\left(
\begin{array}{cc}
\theta & b_2 \\
0 & 0
\end{array}
\right)\in Hom(\CrF,\CrF\otimes\Omega),
$$
with $\theta$ a Higgs field on $\CrE$, and
$b_2\in A^{1,0}(X,E)$ a smooth (1,0)-form. Let $b=b_1+b_2$.
Since $(\CrF,\Theta)$ is a Higgs bundle, we have $\delb_F^2=0$,
$\delb_{Hom(\CrF,\CrF\otimes\Omega)} \Theta=0$
and $\Theta_\wedge \Theta=0$, and this implies that $D''_E b=0$,
where $D''_E=\delb_E+\theta$. Then $b$ defines a class $\overline b\in
H^1_{Dol}(E)$, and under the natural isomorphism we obtain an element
$\xi\in H^1_{DR}(E)$.

\end{proof}

The bijection between 2 and 3 can also be obtained using \cite[lemma
3.5]{S2}. We will finish this section with some remarks, but first 
we need the following lemma

\begin{lemma}
\label{lemma}
If $\theta=0$, then there is a natural isomorphism 
$$
H^i_{Dol}(E) \isom \bigoplus_{j=0}^i H^j(\CrE \otimes \Omega^{i-j})
$$
\end{lemma}

\begin{proof}
Since $\theta=0$,  $D''_E=\delb_E$ and then 

\begin{eqnarray*}
H^i_{Dol}(E) = \frac{\ker(A^i\too A^{i+1})}
{\im(A^{i-1}\too A^i)}= \\
\frac{\oplus \ker(A^{i-j,j} \too A^{i-j,j+1})}
{\oplus \im(A^{i-j,j-1} \too A^{i-j,j})}\isom
\bigoplus_{j=0}^i H^j(\CrE \otimes \Omega^{i-j})
\end{eqnarray*}

\end{proof}

If we consider affine representations in which the linear part is
unitary, then $\theta=0$ and by the previous lemma $H^1_{Dol}(E)=
H^1(\CrE)\oplus H^0(\CrE\otimes \Omega)$. Then such a representation
corresponds to a triple $(\CrE,\xi_1,\xi_0)$ where $\xi_1\in H^1(\CrE)$
and $\xi_0\in H^0(\CrE\otimes \Omega)$. Note that $\xi_1$ is the class
of the extension in item 3 of theorem \ref{affinethm}.

Metrics on holomorphic extensions and the corresponding metric
connections have been studied in \cite{BG} and \cite{DUW}. 
Those connections are different from our connection.
Since they
consider unitary connections, if the extension is not trivial, their
connection never respects the subbundle $E$ (i.e. the image of
$\SA^0(X,E)$ is not in $\SA^1(X,E)$), as opposed to what happens 
with our connection.

Finally let's compare theorem \ref{affinethm} with Simpson's extensions
of Higgs bundles to include arbitrary (not semisimple)
representations \cite{S1}. Take an affine representation, and assume that its
linear part is semisimple. Recall that the 
inclusion $\afn \subset \glp$ gives a $\glp$ representation.
Its semisimple part is just the linear part plus a trivial
one-dimensional factor, so it gives a polystable
Higgs bundle $(\CrE,\theta)\oplus (\CrO_X,0)$. To take into account the
non-semisimple part (the translation part in our case), Simpson adds 
a smooth 1-form $\eta\in A^1(X,\SEnd(E\oplus O_X))$ with 
$D'\eta=0$ and $D''\eta +\eta_\wedge \eta=0$.  
In theorem \ref{affinethm} we construct explicitly this form and show
that it lies in $A^1(X,E)
\subset A^1(X,\SEnd(E\oplus O_X))$, the ``upper-right part'' of the
endomorphisms, and that it is a cocycle, hence defines an element
of the cohomology group.

One can adapt theorem \ref{affinethm} to consider 
parabolic representations of the
fundamental group, i.e. representations $\rho:\pi_1 \to P$ into a 
parabolic subgroup $P$
of $\gln$. If $P=U\cdot L$ is a Levi decomposition of $P$, the
unipotent part $U$ plays the same role as the translation part for
affine representations and the Levi factor $L$ corresponds to the 
``linear part''.
Instead of the short exact sequences considered on item 2, 
one has to consider 
filtrations of flat vector bundles. If the connections induced on the
quotients are semisimple, they give polystable Higgs bundles.  
To recover the parabolic representation from this data, on has
to add a 1-form, analogous to the 1-form defining the de Rham
cohomology in item 4 of the theorem (this 1-form corresponds to the
unipotent part of the representation). Equivalently, one can define
a Higgs field as in item 3, obtaining a semistable (but in general not
polystable) Higgs bundle. Details are given in the appendix.

\section{Appendix: Parabolic representations}
\label{parabolic}

Let $P$ be a parabolic subgroup of $\gln$.
In this appendix we will consider representations of the fundamental
group into $P$. 

Recall that there is a one to one correspondence between parabolic
subgroups of $\gln$ and flags in $\CC^n$ (the parabolic group
associated to a flag is the subgroup of $\gln$ that respects the flag).
Parabolic representations appear naturally when we
consider non-semisimple linear representations.
Let $\rho:\pi_1 \to \gln$ be such a non-semisimple
representation.
Then there is a parabolic subgroup $P$ of $\gln$ such that $\rho$
factors
$$
\rho:\pi_1 \to P \subset \gln,
$$
and such that if the parabolic group $P$ corresponds to a flag 
$$
0=V_0 \subset V_1 \subset V_2 \subset \cdots \subset V_p = \CC^n,
$$
then the associated representation $\rho^{ss}$ in 
$$
\frac{V_1}{V_0} \oplus \frac{V_2}{V_1} \oplus \cdots \oplus 
\frac{V_p}{V_{p-1}}
$$
is semisimple.
Let $P=U\cdot L$ be a Levi decomposition, where $U$ is the maximal
unipotent subgroup of $P$ and $L$ is the Levi factor.
In matrix form, $U$ and $L$ are respectively  matrices of the form
$$
\left(
\begin{array}{cccc}
\id_1 & * & \cdots & * \\
 & \id_2 & \cdots & * \\
 &  & \ddots & \vdots \\
0 &  &  & \id_p 
\end{array}
\right)
\qquad
\left(
\begin{array}{cccc}
*_1 &  &  & 0 \\
 & *_2 &  &  \\
 &  & \ddots &  \\
0 &  &  & *_p 
\end{array}
\right)
$$
where $*_i$ is a  matrix of dimension $\dim (V_i/V_{i-1})$, and
$\id_i$ is the identity matrix of the same dimension.
There is a projection $P\to L$, and the semisimple representation
$\rho^{ss}$ corresponds to the composition $\pi \to P \to L$. 

Since the representation $\rho^{ss}$ is semisimple, it gives a polystable
Higgs bundle $(\CrE_1,\theta_1) \oplus (\CrE_2,\theta_2) \oplus \cdots
\oplus (\CrE_p,\theta_p)$. 
We want to find what extra object we have to give to take into account
the ``non-semisimple part'' $U$.

\subsection{Nonlinear de Rham cohomology}\hfill
\smallskip

Let $(E_i,D_i)$, $i=1,\ldots,p$ be $p$ flat bundles. Let $\uni$
be the subsheaf of $\SEnd(\oplus E_i)$ such that its local sections
respect the filtration 
$$
E_1\subset E_1\oplus E_2 \subset \cdots \subset (E_1\oplus \cdots\oplus
E_p)
$$
and such that they induce the zero endomorphism on each factor $E_i$.
In other words, $\uni$ is the sheaf of Lie algebras of the sheaf of
unipotent groups associated with this filtration.

A flat filtration of vector bundles with quotients $(E_i,D_i)$ 
will be the following data:
\begin{itemize}
\item A filtration of vector bundles
$$
0 = F_0 \subset F_1 \subset F_2 \subset \cdots \subset F_p=F.
$$

\item Isomorphisms $F_i/F_{i-1}\isom E_i$

\item A flat connection $D_F$ on $F$ that respects the filtration and
such that it induces on each quotient $E_i$ (via the previous
isomorphisms)
the connection $D_i$. 

\end{itemize}

If we choose smooth splittings of the
filtration (i.e. a smooth isomorphism $E_1\oplus \cdots \oplus E_p
\isom F$), then we can decompose $D_F = D_L + \eta$ where
$D_L=D_1\oplus \cdots \oplus D_p$ and $\eta\in A^1(X,\uni)$ is a
1-form with values in $\uni$. In matrix form:
$$
D_F=
\left(
\begin{array}{cccc}
D_1& & &0  \\
{} & D_2 \\
 & {} & \ddots \\
0 & & & D_p 
\end{array}
\right)
+
\left(
\begin{array}{cccc}
0 & \eta_{1,2} & \cdots &\eta_{1,p} \\
 & 0 & \cdots & \eta_{2,p} \\
 &   & \ddots & \vdots \\
0&   &    & 0 
\end{array}
\right)
$$
Flatness of $D^{}_F$ translates into $D^{}_L \eta+ \eta_\wedge\eta=0$.
We will think of this as a nonlinear closedness condition.
If we choose a different splitting, the new splitting is isomorphic to
the old one by an isomorphism that has the matrix form
$$
M_0=\left(
\begin{array}{cccc}
\id_1 & M_{1,2} & \cdots &M_{1,p} \\
 & \id_2 & \cdots & M_{2,p} \\
 &   & \ddots & \vdots \\
0&   &    & \id_p 
\end{array}
\right)
$$
and the new 1-form is $\tilde\eta= M^{-1}_0 D^{}_L M^{}_0 + M^{-1}_0 \eta
M^{}_0$. Define the nonlinear de Rham cohomology set as the set of
(nonlinear) closed 1-forms modulo the equivalence relation generated
by the change of splitting
$$
H^1_{DR}(E_\bullet,D_\bullet)=\frac
{\{\eta \in A^1(X,\uni):\quad D^{}_L\eta+\eta_\wedge\eta=0\}}
{\eta \thicksim M^{-1}_0 D^{}_L M^{}_0 + M^{-1}_0 \eta M^{}_0}
$$
Then the flat filtration $F$ gives a well defined element 
 of this
cohomology. Note that this is a pointed set (the distinguished point
being the point corresponding to $\eta=0$). If $p=2$ then it is easy
 to see that this is the (linear) de Rham cohomology $H^1_{DR}(\SHom(E_2,E_1))$
 defined by Simpson \cite{S1}.

There are two different natural definitions of morphism (and
then two different notions of isomorphisms) for
flat filtrations. Let $(F^{}_\bullet,D_F)$ and $(F'_\bullet, D_{F'})$
be two flat filtrations. A \textit{strong} morphism
(resp. isomorphism) of flat filtrations 
is a morphism (resp. isomorphism) $\varphi:F\to F'$ such that  
$\varphi^*D_{F'}=D_F$, (hence $\varphi$ respects the
filtration) and it induces the identity map
on the quotients 
$$
\varphi^{}_i:\frac{F_i}{F_{i-1}}\isom E_i \stackrel{\id}{\too} E_i
\isom\frac{F'_i}{F'_{i-1}}
$$
If the induced maps on quotients $\varphi^{}_i$ are isomorphisms
with $\varphi^*_i D'_i=D^{}_i$ (but $\varphi^{}_i$ is not necessarily
the identity), then
we say that $\varphi$ is a \textit{weak} morphism (resp. isomorphism).

It is easy to check that if $F$ is strongly isomorphic to $F'$, then
the corresponding points in the nonlinear de Rham cohomology
$H^1_{DR}(E_\bullet,D_\bullet)$ 
are the
same $[\eta]=[\eta']$. If they are only weakly isomorphic, then  
if $M$ is the induced map on 1-forms we have $[\eta']=
[M^{-1} D^{}_L M + M^{-1} \eta M]$, but since the induced maps on the
quotients $E_i$ are not necessarily identities, this class is in 
general not equal to $[\eta]$.

\begin{remark}
\textup{
If 
instead of flat connections we use $\delb$-operators, we obtain
the nonlinear $\delb$-cohomology set $H^1_\delb
(E_\bullet,\delb_\bullet)$, and this set parametrizes holomorphic
filtrations with fixed quotients $(E_i,\delb_i)=\CrE_i$. 
If $p=2$ then this filtrations are just extensions of $\CrE_2$ by
$\CrE_1$, we have $\eta_\wedge\eta=0$, and
$H^1_\delb(E_\bullet,\delb_\bullet)=H^1(\SHom(\CrE_2,\CrE_1))$,
the usual cohomology group of coherent sheaves.
}
\end{remark}

\subsection{Main theorem for parabolic representations}\hfill
\smallskip

Let $(\CrE_i,\theta_i)$, $i=1,...,p$ be a collection of polystable 
Higgs bundles with vanishing Chern classes. A filtered Higgs bundle
with quotients $(\CrE_i,\theta_i)$ is a Higgs bundle $(\CrF,\Theta)$
together with a holomorphic filtration 
$$
0 = \CrF_0 \subset \CrF_1 \subset \CrF_2 \subset \cdots \subset 
\CrF_p=\CrF
$$
(such that the Higgs field $\Theta$ respects this filtration), 
and isomorphisms from $(\CrE_i,\theta_i)$ to the induced Higgs 
bundles on the quotients $\CrF_i/\CrF_{i-1}$.
A weak isomorphism of filtered Higgs bundles is an isomorphism of Higgs
bundles, respecting the filtrations, and inducing isomorphisms in the
quotients. If it induces the identity morphisms on the quotients, then
it is called a strong isomorphism.

\begin{theorem}
\label{parathm}
There are natural bijections among the following sets
\begin{enumerate}
\item
The set of $P$ representations such that the induced $L$
representation is semisimple, modulo conjugation by elements of $P$.

\item
The set of weak isomorphism classes
of flat filtrations of vector bundles inducing semisimple flat
connections $D_i$ on the quotients $E_i=F_i/F_{i-1}$
$$
\{(F_\bullet,D_F): D_i\quad \text{semisimple}\}/\isom_{weak}.
$$
\item
The set of weak isomorphism classes of filtered Higgs bundles
inducing polystable Higgs bundles on the quotients $\CrE_i=\CrF_i/
\CrF_{i-1}$
$$
\{(\CrF_\bullet,\Theta): (\CrE_i,\theta_i)\; \text{polystable with
vanishing Chern classes}\}/\isom_{weak}.
$$

\item
The set of isomorphism classes of objects
$$
\{(\CrE_i,\theta_i),\xi\}/\isom
$$
where $(\CrE_i,\theta_i)$ are polystable Higgs bundles, $\xi\in
H^1_{DR}(E_\bullet,D_\bullet)$ (where $(E_i,D_i)$ is the flat bundle
associated to the Higgs bundle $(\CrE_i,\theta_i)$ via a harmonic
metric),
and two such objects are considered
isomorphic if there are isomorphisms $\psi_i:\CrE_i\to \CrE'_i$ of
Higgs bundles sending $\xi$ to $\xi'$.
\end{enumerate}

\end{theorem}

\begin{proof}
$(1 \leftrightarrow 2)$
This is given by holonomy.

$(2 \leftrightarrow 3)$
Follows from \cite[lemma 3.5]{S1} (see the remarks after Simpson's proof).

$(2 \leftrightarrow 4)$ Note that since 
the flat filtration induces semisimple flat connections on the
quotients $E_i$, we get polystable Higgs bundles $(\CrE_i,\theta_i)$,
and as we have already discussed, the flat connection $D_F$ gives a
well defined element of $H^1_{DR}(E_\bullet,D_\bullet)$. It is easy to
check that a weakly isomorphic flat filtration gives isomorphic Higgs
bundles and element $\xi'$.

Conversely, given Higgs bundles $(\CrE_i,\theta_i)$ and $\xi$,
define the filtration $F_\bullet$
$$
0 \subset E_1 \subset E_1\oplus E_2 \subset \cdots \subset 
(E_1\oplus \cdots \oplus E_p).
$$
Take a representative $\eta\in A^1(X,\uni)$ of $\xi$, define
the connection
$$
D_F = (D_1\oplus\cdots \oplus D_p) + \eta,
$$
and this defines a flat filtration.

\end{proof}

\end{document}